\documentclass{amsart}
\usepackage{amsmath,amssymb}
\newtheorem{theorem}{Theorem}
\newtheorem{lemma}[theorem]{Lemma}
\newtheorem{proposition}[theorem]{Proposition}
\newtheorem{corollary}[theorem]{Corollary}

\newtheorem{remark}[theorem]{Remark}

\begin{document}
\title{Commutative Local Rings of bounded module type}
\author{Fran\c cois Couchot}

\begin{abstract} Let $R$ be a local ring of bounded module type. It is
  shown that $R$ is an almost maximal valuation ring if there exists a
  non-maximal prime ideal $J$ such that $R/J$ is an almost maximal
  valuation domain. We deduce from this that $R$ is almost maximal if
  one of the following conditions is satisfied: $R$ is a
  $\mathbb{Q}$-algebra of Krull dimension $\leq 1$ or the maximal
  ideal of $R$ is the union of all non-maximal prime ideals.
\end{abstract}
\maketitle

\bigskip
In this paper, $R$ is an associative and commutative ring with
identity. We will say that $R$ has bounded module type if, for some positive
integer $n$, every finitely generated $R$-module is a direct sum of
submodules generated by at most $n$ elements.
The problem of investigating commutative rings of bounded module
type has been studied by R.B. Warfield \cite{War},  R. Wiegand \cite{Wie},  B.
Midgarden and S. Wiegand \cite{MiWi},   P. Zanardo \cite{Zan},  P. V\'amos \cite{Vam}
and the author \cite{Cou}.
In \cite{War}, R.B. Warfield proved that every local ring of bounded
module type is a valuation ring. By theorems due to D.T. Gill [4]  and
J.P. Lafon \cite{Laf}, a valuation ring is almost maximal if and only if every
finitely generated module is a direct sum of cyclics. So, the following
question can be proposed~:

Is a local ring \ $R$ \ of bounded module type if and only if \ $R$ \
is an almost maximal valuation ring~?

Positive answers are given by  P. Zanardo in \cite{Zan} for the class
of totally branched and discrete valuation domains, by  P. V\' amos in \cite{Vam} for
\ $\mathbb{Q}$-algebra valuation domains and in \cite{Cou} by the author for valuation rings
with a finitely generated maximal ideal.

\bigskip
In this paper, we prove that if \ $R$ \ is a valuation ring of
bounded module type, then \ $R/I$ \ is complete in its ideal topology for every
nonzero and nonarchimedean ideal $I$, and \ $R_J$ \ is almost maximal for every
nonmaximal prime ideal $J$. To obtain these results, as in \cite{Cou}, we adapt to the
nondomain case  Zanardo's methods used in \cite{Zan}. Moreover, we extend results
obtained by  P. V\' amos in \cite{Vam} : every valuation ring \ $R$ \ of bounded
module type is almost henselian, and if \ $R$ \ is a \ $\mathbb{Q}$-algebra with Krull
dimension greater than one, then \ $R$ \ is almost maximal. Finally we obtain
also a positive answer for valuation rings such that the maximal ideal is the
union of all nonmaximal prime ideals.

\bigskip
 For definitions and general facts about valuation rings and their 
modules we refer
to the book by  Fuchs and  Salce \cite{FuSa}. The symbol \ $A\subset B$ \ denotes
that \ $A$ \ is a subset of $B$, possibly \ $A = B$. We recall some definitions
and results which will be used in the sequel.
An \  $R$-module \ $M$ \ is called \textbf{fp-injective} if and only if \
$\mathrm{Ext}^1(F,M) =0$  \ for every finitely presented $R$-module
 $F$. A \textbf{self \ fp-injective ring} \ $R$ \ is \ fp-injective as \
 $R$-module. For an \ $R$-module $M$ and an ideal \ $I$ \ of $R$, we
 set \ $M[I] = \{x\in M \mid I \subset (0 : x)\}$. If \ $M$ \ is \ fp-injective then \ $rM =M[(0:r)]$ \ for every \ $r\in R$.  It is obvious that \ $rM \subset M [(0:r)]$.
Now let
\ $x\in M [(0:r)]$. It follows that we can define an homomorphism \ 
$f : Rr \rightarrow M$
\ by \ $f(tr) = tx$. Since \ $M$ \ is \ fp-injective, \ $f$ \ can 
be extended \
to $R$ \ hence we get \ $x = f(r) =  rf(1)$. Conversely if \ $rM = M[(0 : r)]$,
then every homomorphism \ $f : Rr \rightarrow M$ \ can be extended to $R$, 
hence \ $\mathrm{Ext}^1(R/rR, M) = 0$.
We recall that \ $R$ \ is a \textbf{valuation ring} if and only if its
 ideals are totally ordered by inclusion.

\begin{proposition} \label{P:fp-inj}   Let \ $R$ \ be a valuation ring
 and \ $P$ \ its  maximal ideal. Then :

\begin{enumerate}
\item $R$ \ is self \ fp-injective if and only if each nonunit of \
$R$ \ is a zero-divisor.

In this case, if \ $I$  and $J$ \ are ideals of \ $R$ \ the following 
assertions
hold~:

\item
 If \ $(0 : P) \not= 0$ \ then \ $I = \bigl( 0 : (0 :
I)\bigr)$ \ and if \ $I\varsubsetneq J$ \ we have \ $(0 : J)
\varsubsetneq (0 : I)$.

\item
 If \ $(0 : P) = 0$ \ then \ $I \varsubsetneq \bigl( 0 :
(0 : I)\bigr)$ \ (respectively \ $I\varsubsetneq J$ \ and \ $(0 : J) = (0
: I)$) \ if and only if there exists \ $r\in R$ \ such that \ $I = Pr$ \
and \ $\bigl( 0 : (0 : I)\bigr) = Rr$ \ (respectively $J = Rr$).

\item
 If \ $I\varsubsetneq J$ \ then \ $\bigl( 0 : (0 : I)\bigr)
\subset J$.
\end{enumerate}
\end{proposition}
\textbf{Proof.}   From \cite[Lemma 3]{Gil}  it follows that every
nonunit of \ $R$ \ is a 
zero-divisor if and only if \ $Rr = \bigl( 0 : (0 : r)\bigr)$ \ for 
every \ $r\in R$.
Since \ $R$ \ is a valuation ring then every finitely presented
module is a direct sum of cyclic modules (\cite[Theorem 1]{War} or
\cite[Proposition II.2]{Laf}). 
Hence \ $R$ \ is self
\ fp-injective if and only if \ $\mathrm{Ext}^1 (R/r, R) = 0$ \
for every \ $r\in R$.  From above we deduce that this last condition
is equivalent to \ $\bigl( 0 :(0:r)\bigr) = Rr$ \ for every \ $r\in
R$. 
The assertions 2.  3.  and 4.  follow from  of \cite[Proposition
1.3]{KlLe}.\qed

\bigskip
From now on, $R$ will denote a valuation ring, $P$ its maximal ideal and
\ $E$ \ the \ $R$-injective hull of $R$.

\begin{lemma} \label{L:arch} Let \ $I$\  \  be a nonzero proper ideal
  of \  $R$ \  and \  $p\in P$. \  Then~:
\begin{enumerate}

\item $pI = I$ \ if and only if $(I : p) = I$.

\item $\forall r\in P$, \ $pI\not= I$  and $rI \not= 0$
implies \ $prI \not= rI$.

\item If \ $pI = I$ \ then :
\begin{itemize}
\item[i)] $\forall n\in\mathbb{N}$ \ \ $p^nI = I$ \ and \ $(I : p^n) = I$

\item[ii)] $\forall a \in I$, \ \ $p(Ra : I) = (Ra : I)$

\item[iii)] $\forall a \in I$, \ if \ $p^ma \not= 0$ \  then \
$(Ra : I) = (Rap^m : I)$.
\end{itemize}
\end{enumerate}
\end{lemma}

\textbf{Proof.}

1) Suppose \ $pI = I$.  Let $r\in (I : p)$. If \  $rp = 0$, from
\ $pI
\not= 0$ \ it follows that
$(0 : p) \subset I$ and $r\in I$. If \  $rp \not= 0$ then $rp = tp$ 
for $t\in I$
\ whence we get \  $Rt = Rr$ and $r\in I$.

Suppose $(I : p) = I$.  Then \  $p\notin I$ \  whence for each \ $r\in I$ \
$\exists t\in P$ such that $r = pt$.  We have $t\in (I : p) = I$.

2) If \ $prI = rI$ then $\forall a \in I$, $\exists b\in I$ such
that $rpb = ra$.  If $ra = 0$,  we have $a\in (0 : r) \subset pI$
since $rI \not= 0$.  If $ra \not= 0$,   we
obtain that $Ra = Rpb \subset pI$.

3) For  iii) see the proof of \cite[Lemma 2.3]{Cou}. \
ii) Let
$r\in \bigl( (Ra : I) : p\bigr)$. Then $prI \subset Ra$ and $rI \subset
Ra$. Hence $\bigl( (Ra : I) : p\bigr) = (Ra : I)$  and by using 1) we have
$p (Ra : I) = (Ra : I)$.
\qed

\bigskip
If $I$ is an ideal of $R$,  by using Lemma~\ref{L:arch},  we follow
\cite{FuSa} by calling $I$ \textbf{archimedean ideal} if $(I : p)
\not= I$, for every $p\in P$. Then $I$ is archimedean if and 
only if $R/I$ is a self fp-injective ring (Proposition~\ref{P:fp-inj}).

\bigskip
$R$ is called \textbf{maximal} if and only if every totally ordered
family ${\mathcal{F}}$ of cosets  $\{ r_\lambda + I_\lambda \mid \lambda \in 
\Lambda\}$ has a nonempty intersection, and $R$ is called almost
maximal if the above condition holds whenever  $I = \displaystyle\bigcap_{\lambda\in\Lambda} I_\lambda \not=0$. When  $I = Pr$ \ for some \  $r\in R$, then there exists  $\lambda \in
\Lambda$ such that  $I = I_\lambda$ (else $r \in
\displaystyle\bigcap_{\lambda\in\Lambda} I_\lambda$). In this case we 
deduce that
${\mathcal{F}}$ has a nonempty intersection.

\bigskip
The \textbf{ideal topology} of \  $R$ \ is the linear topology where
the family $\mathcal{J}$  of all nonzero ideals of  $R$  is a base of
neighborhoods  about 0. We said that
$R$ is complete in its ideal topology if and only if the canonical homomorphism
$\phi : R\rightarrow \varprojlim_{I\in\mathcal{J}}R/I$
\ is an isomorphism. If
$\mathcal{F} = \{ I_\lambda \mid \lambda\in\Delta\}$  is a family of 
nonzero ideals
such that  $\displaystyle\bigcap_{\lambda\in\Lambda} I_\lambda = 0$, then it
is easy to verify that  $\mathcal{F}$ is also a base of neighborhoods 
about 0 in the
ideal topology. Consequently  $R$ is complete in its ideal topology 
if and only if
every totally ordered family of cosets  $\{ r_\lambda + I_\lambda \mid \lambda
\in\Lambda\}$,  with $\displaystyle\bigcap_{\lambda\in\Lambda} 
I_\lambda = 0$,
has a nonempty intersection. Thus  $R$  is maximal (respectively 
almost maximal) if
and only if  $R/I$  is complete in its ideal topology for each proper ideal
(respectively nonzero proper ideal) $I$, such that  $I\not= Pr$ for 
any  $r\in R$.

\bigskip
Let now  $e$  be in $E$ not in $R$. As in \cite{SaZa} the breadth
ideal  $B(e)$ of $e$ is defined as follows~: $B(e) = \{ r\in R\mid
e\notin R + rE\}$. The referee suggested me the following proposition that is similar to \cite[Proposition 1.4]{SaZa}.

\begin{proposition} \label{P:breath}  Suppose  $R$ is self fp-injective and let $I$ be a
proper ideal of  $R$, such that  $I \not= Pr$  for every  $r\in R$. 
Then  $R/I$  is
not complete in its ideal topology if and only if there exists  $e\in 
E \setminus
R$ such that   $I = B(e)$. Moreover we have the following properties :
\begin{enumerate}
\item $I = B(e) = (0 : (R : e))$,
\item $(R : e) = P(0 : I)$ and  $(R : e)$  is not finitely generated,
\item and $e$ can be chosen such that  $(0 : e) = 0$.
\end{enumerate}
\end{proposition}

\textbf{Proof.}
Since \  $R/I$ \   is not complete in its ideal topology there exists a
totally ordered family of cosets  $\{ r_\lambda + I_\lambda \mid \lambda \in
\Lambda\}$  with empty intersection such that  $I =
\displaystyle\bigcap_{\lambda\in\Lambda} I_\lambda$. Let $J =
\displaystyle\bigcup_{\lambda\in\Lambda} (0 : I_\lambda)$. We define 
$f : J\rightarrow R$
by  $f(c) = r_\lambda c$ for each  $c$  in  $(0 : I_\lambda )$. If \ 
$f$ \ can be
extended to an endomorphism of  $R$, then, by the same proof (the 
beginning of the
proof) as in \cite[Theorem 2.3]{KlLe} we obtain that
$\displaystyle\bigcap_{\lambda\in\Lambda} r_\lambda + I_\lambda \not= 
\phi$. Thus
$f$ cannot be extended to  $R$, but there exists  $e\in E\setminus R$ such that
$f(c) = ce$ for each  $c$ in $J$. It is obvious that  $J\subset (R : 
e)$. Let $r\in
(R : e)$. We consider the homomorphism
$g : Rr \rightarrow R$ defined by  $g(r) = re$. Since  $R$  is self 
fp-injective, $g$
can be extended to  $R$, hence there exists  $u\in R$ such that 
$g(r) = re = ru$.
If  $J\subset Rr$, then  $g$ is an extension of  $f$  and we obtain a
contradiction. Thus  $Rr \varsubsetneq J, \ J = (R : e)$ and  $(R : e)$  is
not finitely generated. Since  $I \varsubsetneq I_\lambda$ for each 
$\lambda \in
\Lambda$, then  $J\subset (0 : I)$. If  $(0:I)$  is not finitely 
generated then  $J
\subset (0 : I) = P(0 : I)$. If  $(0 : I)$  is finitely generated, then
$J\varsubsetneq (0 : I)$ hence  $J\subset P(0 : I)$. Now the
following equality and inclusions hold : $I = (0 : (0 : I)) \subset (0 : J)
\subset \displaystyle\bigcap_{\lambda\in\Lambda} (0 : (0 : I_\lambda 
))$. Let  $r
\in R\setminus I$. Then there exists $\lambda \in \Lambda$ such that $r\notin
I_\lambda$. If  $I_\mu \varsubsetneq I_\lambda$ for some  $\mu 
\in\Lambda$, then
$(0 : (0 : I_\mu )) \subset I_\lambda$. We get that  $r\notin
\displaystyle\bigcap_{\alpha\in\Lambda} (0 : (0 : I_\alpha))$ and the 
equalities  $I
= (0 : J) = (0 : (R : e))$. We deduce that  $(0 : (0 : J)) = (0 : I)$. When
$J\not=  Pr$  for any   $r\in R$, we have  $J = (0 : I) = P(0 : I)$ 
since  $J$  is
not finitely generated. If \  $J = Pr$\  for some  $r \in R$, we get 
also \  $J =
P(0 : I)$.

Conversely let  \ $e\in E\setminus R$\  and \  $f : (R : e) \rightarrow R$ \  be the
homomorphism defined by  $f(r) = re$. If  $f$  can be extended to 
$g : R\rightarrow R$,
then we obtain that  $(R : e) \subset (0 : e-g(1))$. But  $(R : 
e-g(1)) = (R : e)$
and since  $E$  is an essential extension of  $R$, there exists  $r\in (R : e)$
such that  $r(e-g(1)) \not= 0$.  Thus we get a contradiction. Hence \ 
$f$ \ cannot
be extended to \ $R$ \ and \ $(R : e)$ \ is not finitely generated. Let \ $\{
c_\lambda\mid \lambda \in \Lambda\}$ \ be a set of generators of \ $(R : e)$.
Since \ $R$ \ is self \ fp-injective the restriction of \ $f$ \ to 
\ $Rc_\lambda$
\ can be extended to $R$.  Then there exists \ $r_\lambda\in R$ \ such that \
$f(c_\lambda) = r_\lambda c_\lambda$. Let \ $I_\lambda = (0 : 
c_\lambda)$.  If \
$c_\lambda \in Rc_\mu$, then \ $(r_\lambda - r_\mu) c_\lambda = 0$ \ hence \
$r_\mu \in r_\lambda + I_\lambda$.  Assume that \
$\displaystyle\bigcap_{\lambda\in\Lambda} (r_\lambda + I_\lambda) 
\not= \phi$.  By
the same proof (the end of the proof) as in \cite[Theorem 2.3]{KlLe} we 
deduce that \ $f$
\ can be extended to $R$, hence we get a contradiction. Let \ $I =
\displaystyle\bigcap_{\lambda\in\Lambda} I_\lambda$.  Thus \ $R/I$ \
is not complete in  its ideal topology. The inclusion \ $\bigl( 0 : (R : 
e)\bigr) \subset
I$ \ is obvious. Let \ $r\notin \bigl( 0 : (R : e)\bigr)$.  Then there exists \
$\lambda\in\Lambda$ \ such that \ $rc_\lambda \not= 0$.  Consequently 
\ $r\notin
(0 : c_\lambda)$ \ and we get \ $\bigl( 0 : (R : e)\bigr) = I$.

To complete the proof, we must prove that \ $B(e) = \bigl( 0 : (R : 
e)\bigr)$. Let
 $r\in R\setminus B(e)$.  Then there exist \ $u\in R$ \ and \ $x\in E$ \
such that \ $e = u+rx$.  For each \ $t\in (0:r)$ \ we have \ $te = 
tu$.  If \ $(R :
e) \subset (0 : r)$ \ then the homomorphism \ $f : (R : e) \rightarrow R$ \ 
defined by \
$f(t) = te$ \ can be extended to $R$. It is not possible. Thus \ $(0 :
r)\varsubsetneq (R : e)
$. Since \ $(R : e)$ \ is not finitely generated we get \ $\bigl( 0 : 
(R : e)\bigr)
\varsubsetneq Rr$. Conversely let \ $r\in R \setminus \bigl( 0 : (R : 
e)\bigr)$.
If \ $(0 : P) \not= 0$ \ then \ $(0 : r) \varsubsetneq (R : e)$. 
If \ $(0 : P)
= 0$ \ denote \ $I = \bigl( 0 : (R : e)\bigr)$.   Since \ $I\not= Pr$,  we
deduce that \ $(0 : r)    \varsubsetneq (0 : I)$. If \ $(0 : I)$ \ is 
not finitely
generated then \ $(0 : r) \varsubsetneq P(0 : I) = (R : e)$.  If \ 
$(0 : I) = Rt$
\ for some \ $t\in R$, \ then \ $(R : e) = Pt$.  We have not \ $(0 : r) = Pt$,
else \ $rt \in (0 : P)$ \ and \ $rt \not= 0$.  Hence \ $(0 : r) 
\varsubsetneq (R :
e)$.  Let \ $c\in (R : e) \setminus (0 : r)$. Then the restriction to 
\ $Rc$ \ of
the homomorphism \ $f$ \  defined above can be extended to \ $R$ \ 
since \ $R$ \ is
self \ fp-injective. Thus there exists  \ $u\in R$ \ such that \ 
$tu = te$ \ for
each \ $t\in (0 : r)$.  Since \ $E$ \ is injective, then \ $E [(0 : r)] = rE$ \
hence   there exists \ $x\in E$ \ such that \ $e-u = rx$. We obtain 
that \ $r\notin
B(e)$.
As in \cite[proposition 1.3 ii)]{Cou} we have \ $(0 : e) = 0$ \ or \ $\bigl(0 :
(1-e)\bigr) = 0$. Obviously  \ $(R : e) = (R : 1-e)$ \ and \ 
$B(1-e) = B(e)$.
\qed

\bigskip
Recall that an \ $R$-module \ $M$ \ has \textbf{Goldie dimension} \ $n$ \ if the
injective hull \ $E(M)$ \ is a direct sum of \ $n$ \ indecomposable injective
modules. We denote \ $g(M)$ \ the Goldie dimension of \ $M$ \ and \ 
$\mu (M)$ \ the
minimal number of generators of $M$.

\begin{proposition} \label{P:indecom} Let \ $R$ \ be a valuation ring and \ $I$ \ a
nonarchimedean and nonzero ideal of $R$. If \ $R/I$ \ is not complete 
in its ideal
topology then, for every \ $n\in\mathbb{N}^*$, \ there exists an indecomposable \
$R$-module \ $M$ \ with \ $\mu (M) = n+1$ \ and \ $g(M) = n$.
\end{proposition}

\bigskip
The following two lemmas are needed for the proof of this proposition.

\begin{lemma} \label{L:trans} Suppose \ $R$ \ is self \ fp-injective. Let I be a nonzero
ideal of $R$, $x\in E \setminus R$ \ and \ $a\in R\setminus I$.  If 
\ $ax\in IE$
\ then \ $x\in (I : a)E$.
\end{lemma}
\textbf{Proof.}
Suppose \ $ax = 0$. Let \ $b\in I$, $b\not= 0$. Then there exists \
$d\in (I : a)$ \ such that \ $b = ad$.  Since \ $b\not= 0$ \ then \ $(0 : d)
\subset Ra \subset (0 : x)$. Thus there exists \ $y\in E$ \ such that $x = dy$.

If \ $ax \not= 0$, \ there exist \ c$\in I$ \ and \ $y\in E$ \ such 
that \ $ax =
cy$. Since \ $a\notin I$, \ there exists \ $d\in (I : a)$ \ such that 
\ $c = ad$.
Since \ $ax \not= 0$, \ we have \ $(0 : d) \subset Ra \subset (0 : x-dy)$.
Thus there exists \ $z\in E$ \ such that \ $x-dy = dz$.

Hence \ $x\in dE \subset (I : a)E$. \qed

\begin{lemma} \label{L:dual} Assume \ $R$ \ is self fp-injective and $(0 : P)
\not= 0$.  Let \ $e\in E \setminus R$, $(0 : e) = 0$, \ $I 
= B(e)$, $a\in
I$, $a\not= 0$ \ and \  $J = (Ra : I)$.  We suppose that $I$ is a 
nonarchimedean
ideal.  Let \ $p\in P$ \  such that \ $pI = I$. Then :
\begin{itemize}
\item[i)] $\forall b\in J, \ \exists e_b \in R\setminus P$ \ with \
$(e-e_b) \in (Ra : b) E$.

\item[ii)] $IE = \displaystyle\bigcap_{b\in J}\  (Ra : b) E$.

\item[iii)] Let \ $c, d\in R$ \   with \  $c+de\in IE$.  Then $c\in R
p^m$ and
$d\in Rp^m$ for every  $m\in\mathbb{N}^*$.
\end{itemize}
\end{lemma}

\textbf{Proof.}
For i) and ii) see proof of \cite[Lemma 2.4]{Cou}.

iii) $pI = I \ \Rightarrow \ p\notin I$. \ If \ $c\notin I$
\ or \
$d\notin I$  we show as in  \cite[Lemma 2.4. iii)]{Cou} that \ $Rc = Rd$.  There
exists \ $v\in R\diagdown P$ \ such that \ $d = vc$.  Now suppose $\exists
k\in\mathbb{N}^*$ \ such that \ $c\notin R p^k$. So we have  \ $p^k = uc$ \ for
some \ $u\in P$,  and \ $u(c+de) = uc(1+ve) = p^k (1+ve) \in IE$.  By
Lemma~\ref{L:trans} \ $(1+ve) \in (I : p^k) E = IE$.  We deduce that \ $e\in R + IE$.  By
Proposition~\ref{P:breath} we obtain a contradiction.
\qed

\bigskip
\textbf{Proof of proposition~\ref{P:indecom}.} Let \ $r\in I$, \ $r\not= 0$. \ We
replace \ $R$ \ by \  $R/rP$ \ and assume that \ $R$ \ is self \ fp-injective
and \ $(0 : P) \not= 0$. By Proposition~\ref{P:breath}, there exists \ $e\in 
E\setminus R$ \
such that \ $(0:e) = 0$
\ and \ $(R : e) = P(0 : I)$. Since \ $I$ \ is nonarchimedean, there exists \
$p\in P$ \ such that \ $pI = I$.  From Lemma~\ref{L:arch}, we deduce that \ $p 
(0 : I) = (0
: I) = P(0 : I) = (R : e)$.
Let us fix \ $n\in\mathbb{N}^*$.
Since \ $p^{2(n-1)} I = I$,  there exists $a\in I$ \
such that \ $p^{2(n-1)} a \not= 0$.  For every $k$, \ $1\leq k \leq n$,  \
we define \ $A_k = R ap^{2(k-1)}$. Let us now define $n$ elements of $E$
not in $R$ in the following way : $e_1 = e$  and for every $k$, \ $2 \leq k
\leq n$, \ $e_k = 1 + p^{k-1}e$.  Then we have  $(R : e_k) = (R : p^{k-1}e)
= \bigl( (R : e) : p^{k-1}\bigr)$.  By Lemma~\ref{L:arch}(1 and 3ii) \ 
$(R : e_k) =
(R : e)$  and  $B(e_k) = I$  \ for every $k$, \  $1 \leq k \leq n$.
\ Let \ $J = (Ra : I)$.  Then by lemma~\ref{L:arch}, \ $J = (A_k : I),
\forall k$, $1\leq k\leq n$.

By using Lemma~\ref{L:dual}, for every integer $k$, \ $1\leq k\leq n$,
there exists a family \ $\{ e_k^b \mid b\in J\}$ \ of units of $R$, such that
$(e_k-e_k^b)
\in (A_k : b) E$.

Now we define an $R$-module $M$ with \ $\{ x_0, x_1,\ldots ,
x_n\}$ \ as spanning set and with the following relations~:

-- $(0 : x_k) = A_k$ \quad for every $k$, \ $1\leq k
\leq n$

-- $(0 : x_0) = A_n$

-- $\forall b\in J$, \ \ $b x_0 = b \ \Bigl(
\displaystyle\sum_{k=1}^n \ e_k^b x_k\Bigr)$.

We prove that \  $M$ \  is indecomposable, with \ $\mu (M) =
n+1$ and $g(M) = n$, as in \cite[Proposition 2.6]{Cou},  where
\cite[Lemma 2.4]{Cou} is replaced by Lemma~\ref{L:dual}. \qed

\bigskip
As  P. V\' amos in \cite{Vam}, we call a local ring \ $R$ \textbf{almost
henselian} if every proper factor ring is a henselian ring.

\begin{theorem} \label{T:complet}  Let \ $R$ \ be a local ring of bounded module type.
Then we have the following~:
\begin{enumerate}
\item $R$ is a valuation ring.
\item
\begin{itemize}
\item[a)] For every nonzero and nonarchimedean ideal \ $I$ \ of \ $R$, \
$R/I$ \ is complete in its ideal topology.
\item[b)] If \ $0$ \ is not prime and nonarchimedean then \ $R$ \ is
also complete.
\end{itemize}
\medskip
\item $R$ is an almost henselian ring.
\medskip
\item For every nonmaximal prime ideal $J$, \ $R_J$ \ is almost
maximal.
\end{enumerate}
\end{theorem}

\textbf{Proof.} 1. See \cite[Theorem 2]{War}.

2.a) is a consequence of proposition~\ref{P:indecom}.

2.b) Let \ ${(I_\alpha )}_{\alpha\in\Lambda}$ \ be a family of
ideals such that \ $\displaystyle\bigcap_{\alpha\in\Lambda} I_\alpha = 0$.
If \ $t\notin I_\alpha$ \ \ $\forall\alpha\in\Lambda$, \ then we have 
\ $(0 : t)
= \displaystyle\bigcap_{\alpha\in\Lambda} (I_\alpha : t)$ \ and we can easily
prove that \ $0$ \ nonarchimedean implies \ $(0 : t)$ \ 
nonarchimedean. By
the same proof as in \cite[Proposition 1]{Gil} () we obtain that \ $R$ \ is 
complete if \
$R$ \ is not a domain.

3. Let \ $J$ \ be the nilradical of $R$. From \cite[Theorem 8]{Vam} 
we deduce that \ $R/J$ \ is almost henselian. If \ $J \not= 0$,
then \ $R/J$ \ is complete. Hence \ $R/J$ \ is henselian and since \
$J$ \ is a nilideal, we deduce that \ $R$ \ is henselian.

4. Let \ $I$ \ be a nonzero ideal of \ $R_J, I \subset J
R_J$, \
$\phi : R\rightarrow R_J$ \ the canonical homomorphism and \ $I^{c} = \phi^{-1} (I)$.
If \ $s\in P\setminus J$ \ and \ $r \in I^{c}$, then there exists \ 
$t\in P$ \ such
that \ $r = st$ \ and \ $\frac{t}{1} = \frac{r}{s} \in
I$.  We deduce that \ $sI^{c} = I^{c}$ \ hence \ $I^{c}$ \ is a nonarchimedean
ideal of \ $R$.  Since \ $R/I^{c}$ \ is complete in its ideal topology, by the
same proof as in \cite[Lemma 2]{Gil}  we prove that \ $R_J/I$ \ is also complete.
Hence \ $R_J$ \ is almost maximal.
\qed

\bigskip
\begin{proposition} \label{P:complet} 
Let \ $R$ \ be a valuation ring. Suppose there
  exists a nonmaximal prime ideal \ $J$ \ such that \ $R/J$ \ is
  almost maximal.  Then, for every archimedean ideal $I$, $R/I$ is
  complete in its ideal topology. 
\end{proposition}
\textbf{Proof.}
Suppose there exists an archimedean ideal \ $I$, \
$I\not= Pr$ \ for each \ $r\in R$, \ such that \ $R/I$ \ is not complete in its
ideal topology. If \ $I\not= 0$ \ we can replace \ $R$ \ by \ $R/I$ \ and
assume that \ $I = 0$, \ $R$ \ self fp-injective and \ $(0 : P) 
= 0$.  By
proposition~\ref{P:breath} there exists \ $e\in E\setminus R$ \ such that \ $B(e) 
= 0$. Let
\ $J$ \ be a nonmaximal prime ideal and \ $t\in  P\setminus J$. Then \ $0
\not= (0 : t) \subset J$. Let \ $s\in (0 : t)$, $s\not= 0$. Since 
\ $B (e) =
0$ \ there exist \ $u\in R$ \ and \ $x\in E$ \ such that \ $e = u+sx$. We
deduce easily that \ $(R : x) = s(R : e) = sP$ \ and \ $B(x) = (0 : 
sP) = (0 : s)$.
By proposition~\ref{P:breath} \ $R/B(x)$ \ is not complete in its ideal topology and  \ $J
\varsubsetneq Rt \subset B(x)$. We obtain that \ $R/J$ \ is not almost
maximal. \qed

\bigskip
From this proposition and Theorem~\ref{T:complet}  we deduce the following corollary~:

\begin{corollary} \label{C:almo} Let \ $R$ \ be a local ring of bounded module type. Suppose
there exists a  nonmaximal prime ideal \ $J$ \ such that \ 
$R/J$ \ is almost maximal.

Then \ $R$ \ is an almost maximal valuation ring.
\end{corollary}

\bigskip
Now we can give a positive answer to our question for some classes of
valuation rings.

\begin{theorem} \label{T:main} Let \ $R$ \ be a local ring of bounded module type.
Suppose that one of the three following conditions is verified~:
\begin{enumerate}
\item $P$\  is finitely generated.
\item $R$\  is a \ $\mathbb{Q}$-algebra with Krull dimension $\geq
  1$. 
\item$P$\  is the union of all nonmaximal prime ideals of $R$.
\end{enumerate}
Then \ $R$ \ is an almost maximal valuation ring.
\end{theorem}
\textbf{Proof.}
1) It is the main result of \cite{Cou}. We can deduce easily this result
from Theorem~\ref{T:complet} and Corollary~\ref{C:almo} because, if \ $R$ \ is not artinian, then
\ $J =
\displaystyle\bigcap_{n\in\mathbb{N}} P^n$ \ is a prime ideal and \ $R/J$ \ 
is a discrete
rank one valuation domain. Consequently \ $R/J$ \ is almost maximal.

2) Let \ $J$ \ be a nonmaximal prime ideal of $R$. By \cite[Theorem
8]{Vam} 
$R/J$ \ is almost maximal.

3) Let \ $J$ \ be a nonmaximal
and nonzero prime ideal.  We replace \ $R$
\ by \ $R/J$ \ and assume that \ $R_Q$ \ is maximal for every
nonmaximal prime ideal $Q$.  Let$K$ be the field of fractions of $R$, and $X =\mathrm{Spec} R\setminus \{ P\}$.  If \ $x\in K\setminus R$ \ then \ $x = \frac{1}{s}$ \ where \ $s\in P$.  Since \ $P = \displaystyle\bigcup_{Q\in X} Q$, \
$\exists Q\in X$ \ such that \ $s\in Q$.  We deduce that \ $x\notin 
R_Q$ \ and \
$R = \displaystyle\bigcap_{Q\in X} R_Q$.  By \cite[proposition 4]{Zel} $R$ \ is
linearly compact in the inverse limit topology.  Since \ $R$ \ is a 
Hansdorff space
in this linear topology then every nonzero ideal is open and also 
closed. Hence \
$R$ \ is linearly compact in the discrete topology.
\qed

\bigskip
\begin{remark} \textnormal{Let \ $J = \displaystyle\bigcup_{Q\in X}Q$.  If \
$J\not= P$, \ then \ $R/J$ \ is an archimedean valuation domain. From corollary~\ref{C:almo}, if \ $R$ \ is of bounded module type, then \ $R$ \ is almost maximal if and
only if \ $R/J$ \ is almost maximal. So, to give a definitive answer to our
question, we must solve this problem for archimedean valuation rings.}
\end{remark}

\end{document}